\documentclass{amsart}

\usepackage{amsmath}
\usepackage{amscd}
\usepackage{amssymb} 

\newcommand{\cal}{\mathcal}
\newcommand{\bk}{{\bf k}}

\DeclareMathOperator{\Img}{Im}

\DeclareMathOperator{\Ker}{Ker}

\newtheorem{theorem}{Theorem}[section]
\newtheorem{theorem/definition}{Theorem/Definition}[section]
\newtheorem{proposition}{Proposition}[section]
\newtheorem{lemma}{Lemma}[section]

\theoremstyle{remark}
\newtheorem{remark}{Remark}[section]

\theoremstyle{definition}

\begin{document}

\title[Frobenius Manifold structure on Dolbeault Cohomology]
{Frobenius Manifold Structure on Dolbeault Cohomology 
and Mirror Symmetry}
\author{Huai-Dong Cao \& Jian Zhou}
\address{Department of Mathematics\\
Texas A \& M University\\
College Station, TX 77843}
\email{cao@math.tamu.edu,  zhou@math.tamu.edu}
\begin{abstract}
We construct a differential 
Gerstenhaber-Batalin-Vilkovisky algebra from Dolbeault complex of any closed
K\"{a}hler manifold, 
and a Frobenius manifold structure on  Dolbeault cohomology.
\end{abstract}
\maketitle
\date{}

\footnotetext[1]{Both authors were supported in part by NSF }

String theory leads to the mysterious Mirror Conjecture, see Yau \cite{Yau} 
for the history.
One of the mathematical predictions made by physicists
based on this conjecture is the formula due to 
Candelas-de la Ossa-Green-Parkes \cite{Can-Oss-Gre-Park} 
on the number of rational curves of any degree 
on a quintic in ${\Bbb C}{\Bbb P}^4$. 
Recently, it has been proved by Lian-Liu-Yau \cite{Lia-Liu-Yau}.

The  theory of quantum cohomology, also suggested by physicists,
has lead to a better mathematical formulation of
the Mirror Conjecture.
As explained in Witten \cite{Wit},
there are two topological conformal field theories on 
a Calabi-Yau manifold $X$: 
A theory is independent of the complex structure of $X$,
but depends on the K\"{a}her form on $X$, 
while B theory is independent of the K\"{a}hler form of $X$,
but depends on the complex structure of $X$.
Vafa \cite{Vaf} explained how two quantum cohomology 
rings ${\cal R}$ and ${\cal R}'$ arise from these theories,
and the notion  of mirror symmetry could be translated into the 
equivalence of A theory on a Calabi-Yau manifold $X$
with B theory  on another Calabi-Yau manifold $\hat{X}$,
called the mirror of $X$, such that the quantum ring
${\cal R}$ can be identified with ${\cal R}'$.

For a mathematical exposition  in terms of variation of 
Hodge structures, see e.g. Morrison \cite{Mor}
or Bertin-Peters \cite{Ber-Pet}.
There are two natural Frobenius algebras on any 
Calabi-Yau $n$-fold,
\begin{eqnarray*}
A(X) = \oplus_{k=0}^n H^k(X, \Omega^k), &
B(X) =  \oplus_{k=0}^n H^k(X, \Omega^{-k}),
\end{eqnarray*}
where  $\Omega^{-k}$ is the sheaf of holomorphic sections
to $\Lambda^kTX$.
By Hodge theory, $B(X)$ can be identified with $H^n(X, {\Bbb C})$.
By Bogomolov-Tian-Todorov theorem,
there is a deformation of the complex structures on $X$
parameterized by an open set in $H^1(X, \Omega^{-1})$. 
Therefore, one gets a family of Frobenius algebra 
structures on $H^n(X, {\Bbb C})$. 
Every Frobenius algebra structure can be characterized by a cubic
polynomial $\Phi$
(in Physics literature, it is called the Yukawa coupling),
so we get a family $\Phi_B(X)$ of cubic polynomials on
$H^n(X, {\Bbb C})$
parameterized by an open set in $H^1(X, \Omega^{-1})$.
An additional structure provided by algebraic geometry is 
the Gauss-Manin connection on this family.
It is a flat connection with some extra properties.
On the other hand, counting of rational curves provides
a family $\Psi_A(X)$ of cubic polynomials on $A(X)$
parameterized by an open set in $H^1(X, \Omega^1)$.
One version of Mirror Conjecture 
is the conjectural existence, for a Calabi-Yau $3$-fold $X$,
of another Calabi-Yau $3$-fold $\hat{X}$,
such that one can identify 
$B(X)$ with $A(\hat{X})$,
 $\Phi_B(X) = \Psi_A(\hat{X})$,
and vice versa.

Witten \cite{Wit1} proposed that the associativity property of 
quantum cohomology can be encoded in WDVV 
(Witten-Dijkgraaf-Verlinde-Verlinde) equations (p.275, (3.14)).
Dubrovin \cite{Dub1} studied family of Frobenius algebras
which are parameterized by the underlying vector space,
such that the family of structure cubic polynomials 
admits some symmetry, so that it becomes the third derivative
of a single potential function.
This potential function satisfies WDVV equations.
Such a family will be called a {\em potential family}.
He then introduced and studied the concept of
Frobenius manifolds.
See Dubrovin \cite{Dub2}. 
Ruan-Tian \cite{Rua-Tia} gave a mathematical 
formulation of quantum cohomology (see also Liu \cite{Liu}
and McDuff-Salamon \cite{McD-Sal}).
They also proved  the potential property (WDVV equations) of the
quantum cohomology,
and gave the construction of a flat connection 
(suggested  by Dubrovin's construction), 
which is the candidate for the mirror of 
Gauss-Manin connection.
(In  general, there is some convergence problem 
involved in their construction, except for the case of
complete intersection Calabi-Yau manifolds. See Tian \cite{Tia}, $\S 10$.)
Witten \cite{Wit} suggested two kinds of extended moduli spaces,
one containing the deformation space of the complex structure, 
the other containing the complexified K\"{a}hler cone.
Kontsevich \cite{Kon},  besides reformulating Mirror Conjecture
in a more general setting,
proposed a construction of extended moduli space of
complex structures.
Recently, Barannikov-Kontsevich \cite{Bar-Kon}
constructed an extended formal moduli space for 
any Calabi-Yau manifold, as a generalization of 
Bogomolov-Tian-Todorov theorem.
Furthermore, they showed that there is a structure of formal
Frobenius manifold on it.
Manin \cite{Man} has formulated a generalization 
to the construction of a formal Frobenius manifold from
any differential Gerstenhaber-Batalin-Vilkovisky (dGBV) algebra
with some mild conditions.

By Hodge theory (see e.g. Griffiths-Harris \cite{Gri-Har}),
there is a decomposition
$$H^*(X, {\Bbb C}) \cong \oplus_{p, q}  H^{p, q}(X),$$
where $H^{p, q}(X)$ is the Dolbeault cohomology.
Furthermore, exterior  product induces 
$$H^{p, q}(X) \wedge H^{r, s}(X) \subset H^{p+r, q+ s}(X).$$
In Cao-Zhou \cite{Cao-Zho1}, the authors developed a quantum de Rham
cohomology theory.
In that theory,  quantum Dolbeault cohomology can be defined
such that it satisfies the usual properties of the ordinary 
Dolbeault cohomology.
This suggests the construction of Frobenius manifold structure on
Dolbeault cohomology. 
In this paper, we construct dGBV algebra 
for a general closed K\"{a}hler manifold.
From this, we obtain a formal Frobenius
manifold on the $A$-side for any closed K\"{a}hler manifold.
Comparison with that of Barannikov-Kontsevich \cite{Bar-Kon} suggests
that these two kinds of formal Frobenius manifolds should be
mirror image of each other.
A standard argument in Kodaira-Spencer-Kuranishi deformation
theory shows that one actually obtains a  Frobenius (super)manifold
structure on a neighborhood of the origin in the relevant
cohomology, in the setting of both Barannikov-Kontsevich and ours.
Now we have two formal Frobenius manifolds on the
$A$-side: one from counting rational curves by Gromov-Witten
invariants, 
the other from our construction.
We conjecture that  they can be identified. 
If this is true, it should have some applications 
in enumerative geometry.

{\bf Acknowledgements}. 
{\em The authors would like to thank Gang Tian and S.-T. Yau for their interest.
The work in this paper is carried out
while the second author is visiting Texas A$\&$M University.
He likes to express his appreciation for the hospitality
and financial support of the Mathematics Department 
and the Geometry-Analysis-Topology group.
He also likes to thank Blaine  Lawson for bringing his attention to 
the reference \cite{Kon2}, 
which made him come across Barannikov-Kontsevich \cite{Bar-Kon}
and Manin \cite{Man} on the internet.}

\section{A construction of Frobenius (super)manifolds}

In this section,  we review a construction of Frobenius supermanifolds.
For details, the reader should consult the papers
by Tian \cite{Tia}, Todorov \cite{Tod}, 
Barannikov-Kontsevich \cite{Bar-Kon} and Manin \cite{Man}.
Here, we follow the formulation by Manin \cite{Man}.

The input of the construction is a differential 
Gerstenhaber-Batalin-Vilkovisky
(dGBV) algebra (see e.g. \cite{Man}, $\S 5$). 
Let $({\cal A}, \wedge)$ be a supercommutative algebra 
with identity over a field $\bk$, 
i.e., 
$({\cal A}, \wedge)$ is an algebra with identity over $\bk$, 
furthermore, ${\cal A}$ has ${\Bbb Z}_2$-grading, such that
for any homogeneous elements $a, b \in {\cal A}$ 
with degrees $|a|$, $|b|$ respectively, 
we have
$$a \wedge b = (-1)^{|a|\cdot |b|} b \wedge a.$$
Assume that there are two $\bk$-linear maps of 
odd degrees $\delta$ and $\Delta$ on ${\cal A}$, such that 
\begin{itemize}
\item $\delta^2 = 0$, $\delta\Delta+\Delta\delta = 0$, 
$\Delta^2 = 0$.
\item $\delta$ is a derivation: 
$\delta (a\wedge b) = (\delta a) \wedge b + (-1)^{|a|}a \wedge (\delta b)$.
\item If $[a \bullet b] = (-1)^{|a|}(\Delta(a\wedge b) - 
\Delta a \wedge b - (-1)^{|a|} a \wedge \Delta b)$, then 
$[a \bullet \cdot]: {\cal A} \rightarrow {\cal A}$ is a derivation
of degree $|a| + 1$, i.e.,
\begin{eqnarray*}
[a \bullet (b \wedge c)] = [a \bullet b] \wedge c 
+ (-1)^{(|a|+1)|b|}b \wedge [a \bullet c].
\end{eqnarray*}
\end{itemize}
The  quintuple 
$({\cal A}, \wedge, \delta, \Delta, [\cdot \bullet \cdot])$
with the above properties is  called a dGBV algebra.
As in Manin \cite{Man}, $\S 5$, we have
\begin{eqnarray*}
& [a \bullet b ] = -(-1)^{(|a|+1)(|b|+1)}[b \bullet a], \\
& [a \bullet [b \bullet c]] = [[a \bullet b] \bullet c]+ 
(-1)^{(|a|+1)(|b|+1)} [b \bullet [a \bullet c]],\\
& \Delta [a \bullet b] = [\Delta a \bullet b] + 
(-1)^{|a|+1}[a \bullet \Delta b], \\
& \delta [a \bullet b] = [\delta a \bullet b ] + (-1)^{|a|+1}[a \bullet \delta b].
\end{eqnarray*}
We will be interested in the cohomology ring $H = H({\cal A}, \delta)$,
which we assume  to be finite dimensional.
An {\em integral} on ${\cal A}$ is an even linear functional 
$\int: {\cal A} \rightarrow \bk$, such that
\begin{eqnarray}
& \int (\delta a) \wedge b = 
	(-1)^{|a|+1} \int a \wedge (\delta b), \label{eqn:int1} \\
& \int (\Delta a) \wedge b = 
	(-1)^{|a|} \int a \wedge (\Delta b), \label{eqn:int2}
\end{eqnarray} 
for any homogeneous $a, b \in {\cal A}$.
It follows from $(\ref{eqn:int1})$ that $\int$ induces a
well-defined $\bk$-bilinear functional
\begin{eqnarray*}
& (\cdot, \cdot): 
H({\cal A}, \delta) \otimes H({\cal A}, \delta) \rightarrow \bk, \\
& (\alpha, \beta) = \int \alpha \wedge \beta,
\end{eqnarray*}
where $[\alpha]$ and $[\beta]$ are cohomology classes represented
by $\delta$-closed $a$ and $b$ respectively.
If $(\cdot, \cdot)$ is nondegenerate, 
we say that the integral is {\em nice}.
It is obvious that 
$$(\alpha \wedge \beta, \gamma) = (\alpha, \beta \wedge \gamma).$$
By definition, $(H({\cal A}, d), \wedge, (\cdot, \cdot))$ is then
a {\em Frobenius algebra}, when ${\cal A}$ has a nice integral.

To find a deformation of the ring structure on $H$,
consider $\delta_a: {\cal A} \rightarrow {\cal A}$ for even 
$a \in {\cal A}$.
When $a$ satisfies
\begin{equation} \label{eqn:MC}
\begin{split}
&\delta a + \frac{1}{2}[a \bullet a] = 0, \\
&\Delta a = 0,
\end{split} 
\end{equation}
then $({\cal A}, \wedge, \delta_a, \Delta, [\cdot \bullet \cdot])$ 
is also a dGBV algebra.
If $\int$ is an integral for $({\cal A}, \wedge, \delta, \Delta,
[\cdot \bullet \cdot])$, 
so is it for 
$({\cal A}, \wedge, \delta_a, \Delta, [\cdot \bullet \cdot])$. 
If there is a natural way to identify $H({\cal A}, \delta_a)$ 
with $H = H({\cal A}, \delta)$ which preserves $(\cdot, \cdot)$,  
then we get another Frobenius algebra structure on $H$.
The construction of Frobenius manifold structure 
is based on the existence of a solution 
$\Gamma = \sum \Gamma_n$ to (\ref{eqn:MC}),
such that $\Gamma_0 = 0$,
$\Gamma_1 = \sum x^j e_j$, $e_j \in \Ker \delta \cap \Ker \Delta$.
For $n > 1$, $\Gamma_n \in \Img \Delta$ is 
a homogeneous super polynomial of degree $n$ in $x^j$'s,
such that the total degree of $\Gamma_n$ is even. 
Furthermore, $x^0$ only appears in $\Gamma_1$.
Such a solution is called a {\em normalized universal solution}.
Its existence can be established inductively. 
This is how Tian \cite{Tia1} and Todorov \cite{Tod} proved
that the deformation of complex structures on a Calabi-Yau
manifold is unobstructed.
This was generalized by Barannikov-Kontsevich \cite{Bar-Kon}
to the case of extended moduli space of  complex structures
of a Calabi-Yau manifold.
Manin \cite{Man} further generalized it to the case of 
dGBV algebras, which can be stated as follows:

\begin{theorem} \label{thm:construction}
Let $({\cal A}, \wedge, \delta, \Delta, [\cdot\bullet\cdot])$
be a dGBV algebra which satisfies the following conditions:
\begin{enumerate}
\item $H=H({\cal A}, \delta)$ is finite dimensional.
\item There is  a nice integral on ${\cal A}$. 
\item The inclusions $i: (\Ker \Delta, \delta) \hookrightarrow
({\cal A}, \delta)$ and $j: (\Ker \delta, \Delta) 
\hookrightarrow ({\cal A}, \Delta)$
induce isomorphisms on cohomology.
\end{enumerate}
Then there is a structure of formal Frobenius  manifold 
on the formal spectrum of $\bk[[H']]$, 
the algebra of formal power series generated by $H'$, 
where $H'$ is the dual $\bk$-vector space of $H$.
\end{theorem}

\section{Frobenius manifold structures 
on Dolbeault cohomology}

In this section, $(X, g, J)$ will be a closed 
K\"{a}hler manifold with K\"{a}hler form $\omega$.
Consider the quadruple 
$(\Omega^{*, *}(X), \wedge, \delta=\bar{\partial}, 
\Delta=\partial^*)$.
It is well-known that $\bar{\partial}^2 = 0$,
$(\partial^*)^2 = 0$, and 
$\bar{\partial}\partial^* + \partial^*\bar{\partial} = 0$.
Furthermore,  $\bar{\partial}$ is a derivation.
Set 
$$[a \bullet b]_{\partial^*} = (-1)^{|a|}(\partial^*(a \wedge b ) 
- \partial^* a \wedge b - (-1)^{|a|}a \wedge \partial b).$$
(To the authors' knowledge,
a formula of this type in differential geometry was first discovered 
by Tian \cite{Tia1} (see also Todorov \cite{Tod}) to prove the important
result that deformations of Calabi-Yau manifolds are unobstructed. )
Since both $\omega$ and $J$
are parallel with respect to the Levi-Civita connection,
it follows that near  each $x \in X$, one can find 
a local frame $\{e^1, \cdots, e^n \}$ of $T^{1, 0}X$, 
such that 
$\omega = e^1 \wedge \bar{e}^1 + \cdots + e^n \wedge \bar{e}^n$.
Furthermore, $\nabla_{e_j} e_k = \nabla_{\bar{e}_j}e_k 
= \nabla_{e_j} \bar{e}_k = 
\nabla_{\bar{e}_j} \bar{e}_k = 0$ at $x$.
Let $\{e^1, \cdots, e^n\}$ be the dual basis.
Then at $x$,
for any $\alpha \in \Omega^{*, *}(X)$, we have
(Griffiths-Harris \cite{Gri-Har}, p. 113):
\begin{eqnarray*}
 \partial \alpha = e^j \wedge \nabla_{e_j} \alpha, 
&& \bar{\partial} \alpha 
	= \bar{e}^j \wedge \nabla_{\bar{e}_j} \alpha, \\
 \partial^* \alpha = -e_j \vdash \nabla_{\bar{e}_j} \alpha, 
&& \bar{\partial}^* \alpha 
	= -\bar{e}_j \vdash \nabla_{e_j} \alpha,
\end{eqnarray*}
 
\begin{lemma} $[a \bullet (b \wedge c]]_{\partial^*}
= [a \bullet b]_{\partial^*} \wedge c 
+ (-1)^{(|a|+1)|b|} b \wedge [a \bullet c]_{\partial^*}.$
\end{lemma}

\begin{proof} The left hand side is given by
\begin{eqnarray*}
&   & (-1)^{|a|} (\partial^*(a \wedge b \wedge c)
	-  (\partial^*a) \wedge b \wedge c
	- (-1)^{|a|} a \wedge \partial^* (b \wedge c)) \\
& = & (-1)^{|a|} (- e_j \vdash \nabla_{\bar{e}_j} (a \wedge b \wedge c) 
+ (e_j \vdash \nabla_{\bar{e}_j} a) \wedge b \wedge c \\
&   & + (-1)^{|a|} a \wedge e_j \vdash \nabla_{\bar{e}_j} (b \wedge c)) \\
& = & (-1)^{|a|} (- e_j \vdash
( \nabla_{\bar{e}_j} a \wedge (b \wedge c)
+ a \wedge \nabla_{\bar{e}_j} (b \wedge c))
+ (e_j \vdash \nabla_{\bar{e}_j} a) \wedge b \wedge c  \\
&   & + (-1)^{|a|} a \wedge e_j \vdash \nabla_{\bar{e}_j} (b \wedge c)) \\
& = & (-1)^{|a|} (- (e_j \vdash \nabla_{\bar{e}_j} a) \wedge (b \wedge c)
 - (-1)^{|a|} \nabla_{\bar{e}_j} a \wedge e_j \vdash (b \wedge c) \\
&   &  -(e_j \vdash a) \wedge \nabla_{\bar{e}_j} (b \wedge c)
 - (-1)^{|a|} a \wedge e_j \vdash  \nabla_{\bar{e}_j} (b \wedge c) \\
&  &  + (e_j \vdash \nabla_{\bar{e}_j} a) \wedge b \wedge c
+ (-1)^{|a|} a \wedge e_j \vdash \nabla_{\bar{e}_j} (b \wedge c)) \\
& = & - \nabla_{\bar{e}_j} a \wedge e_j \vdash (b \wedge c)
- (-1)^{|a|} (e_j \vdash a) \wedge \nabla_{\bar{e}_j} (b \wedge c).
\end{eqnarray*}
The right hand side is computed in a similar fashion.
\begin{eqnarray*}
&   & [a \bullet b] \wedge c 
	+ (-1)^{(|a|+1)|b|} b \wedge [a \bullet c] \\
& = & (-1)^{|a|}(\partial^*(a \wedge b) - (\partial^*a) \wedge b
	- (-1)^{|a|} a \wedge \partial^*b) \wedge c \\
&   & + (-1)^{(|a|+1)|b|}(-1)^{|a|} b \wedge 
	(\partial^*(a \wedge c) - (\partial^*a) \wedge c
	- (-1)^{|a|} a \wedge \partial^*c) \\
& = & (-1)^{|a|} (- (e_j \vdash \nabla_{\bar{e}_j} (a \wedge b)) \wedge c 
+ (e_j \vdash \nabla_{\bar{e}_j} a) \wedge b \wedge c \\
&   & + (-1)^{|a|} a \wedge( e_j \vdash \nabla_{\bar{e}_j} b) \wedge c) \\
&   & + (-1)^{(|a|+1)|b|}(-1)^{|a|} b \wedge 
	(-e_j \vdash \nabla_{\bar{e}_j} (a \wedge c) 
	+ (e_j \vdash \nabla_{\bar{e}_j} a) \wedge c  \\
&   &	+ (-1)^{|a|} a \wedge (e_j \vdash \nabla_{\bar{e}_j} c))  \\
& = & -(-1)^{|a|} (e_j \vdash \nabla_{\bar{e}_j} (a \wedge b)) \wedge c 
- (-1)^{(|a|+1)|b|}(-1)^{|a|} b \wedge 
	e_j \vdash \nabla_{\bar{e}_j} (a \wedge c) \\
&   & + a \wedge( e_j \vdash \nabla_{\bar{e}_j} b) \wedge c
+ (-1)^{|b|} a \wedge b \wedge (e_j \vdash \nabla_{\bar{e}_j} c) \\
&   &  + 2 (-1)^{|a|} (e_j \vdash \nabla_{\bar{e}_j} a) \wedge b \wedge c \\
& = & - \nabla_{\bar{e}_j} a \wedge (e_j \vdash b) \wedge c
 - (-1)^{|b|} \nabla_{\bar{e}_j} a \wedge b \wedge (e_j \vdash c) \\
&   & - (-1)^{|a|} (e_j \vdash a) \wedge \nabla_{\bar{e}_j} b \wedge c
- (-1)^{|a|} (e_j \vdash a) \wedge b \wedge \nabla_{\bar{e}_j} c \\
& = & - \nabla_{\bar{e}_j} a \wedge e_j \vdash (b \wedge c)
- (-1)^{|a|} (e_j \vdash a) \wedge \nabla_{\bar{e}_j} (b \wedge c).
\end{eqnarray*}
\end{proof}

Thus $(\Omega^{*, *}(X), \wedge, \delta=\bar{\partial}, 
\Delta=\partial^*)$ is a dGBV algebra.

\begin{remark}
One can obtain an indirect proof as follows.
Cao-Zhou \cite{Cao-Zho2} showed that for a 
K\"{a}hler manifold,
$\sqrt{-1}(\bar{\partial}^* - \partial^*)$
is the operator $\Delta$ defined in Koszul \cite{Kos},
who also defined  the {\em covariant Schouten-Nijenhuis 
bracket} 
$$\{ a, b\} = (\Delta a) \wedge b + (-1)^{|a|} a \wedge (\Delta b) 
	- \Delta ( a \wedge b)$$
which satisfies the property
$$\{a, b \wedge c\} = \{a, b\} \wedge c
	+(-1)^{(|a|+1)|b|} b \wedge \{a, c\}.$$
As pointed out in Cao-Zhou \cite{Cao-Zho2},
$$[a \bullet b ] = (-1)^{|a|}\{a, b \}$$
is the bracket generated by $\Delta$.
Using the type decomposition and 
$[a \bullet (b \wedge c)] = [a \bullet b] \wedge c 
+ (-1)^{(|a|+1)|b|} b \wedge [a \bullet c]$,
one then gets
$$[a \bullet (b \wedge c]]_{\partial^*}
= [a \bullet b]_{\partial^*} \wedge c 
+ (-1)^{(|a|+1)|b|} b \wedge [a \bullet c]_{\partial^*}.$$
\end{remark}

Let $\int_X: \Omega^{*, *}(X) \rightarrow {\Bbb C}$
be the ordinary integration of differential forms.
Analogous to Claim 4.1 in Barannikov-Kontsevich \cite{Bar-Kon}, we have

\begin{lemma} 
The following identity holds on closed K\"{a}hler manifolds:
\begin{eqnarray*}
& \int_X \bar{\partial} a \wedge b 
= (-1)^{|a|+1} \int_X a \wedge \bar{\partial} b, \\
& \int_X {\partial}^* a \wedge b 
= (-1)^{|a|} \int_X a \wedge {\partial}^* b,
\end{eqnarray*}
\end{lemma}

\begin{proof}
By linearity, we can assume that $a$ has type $(p, q)$.
To prove the first identity,
we only need to consider the case of $b$ having type 
$(n-p, n-q - 1)$. 
Then we have $\partial(a \wedge b) = 0$, 
since it has type $(n+1, n-1)$.
Hence
$d(a \wedge b) = \bar{\partial}(a \wedge b)$.
Similarly, we have $d a  \wedge b = \bar{\partial}a \wedge b$.
By Stokes theorem, we have
\begin{eqnarray*}
\int_X \bar{\partial} a \wedge b  
+ (-1)^{|a|} a \wedge  \bar{\partial} b
=  \int_X  \bar{\partial} (a \wedge b) 
= \int_X d (a \wedge b) = 0.
\end{eqnarray*}
To prove the second equality,  
recall that $\partial^* = -*\partial*$, and 
$*^2 \alpha = (-1)^{|\alpha|}\alpha$ 
(see Griffiths-Harris \cite{Gri-Har}, p. 82).
Then since $|b| = |a| -1$, we have
\begin{eqnarray*}
&   & \int_X {\partial}^* a \wedge b 
= (-1)^{|b|} \int_X {\partial}^* a \wedge * (*b) \\
& = & (-1)^{|b|}({\partial}^* a, *b) =(-1)^{|b|} (a, \partial*b) \\
& = &  (-1)^{|b|} \int_X a \wedge *\partial *b 
= -(-1)^{|b|} \int_X a \wedge \partial^*b \\
& = & (-1)^{|a|}  \int_X a \wedge \partial^*b.
\end{eqnarray*}
\end{proof}

\begin{proposition}
The two natural inclusions 
$i: (\Ker \partial^*, \bar{\partial}) \rightarrow 
	(\Omega^{*, *}(X), \bar{\partial})$
and 
$j: (\Ker \bar{\partial},  \partial^*) \rightarrow 
	(\Omega^{*, *}(X), \partial^*)$
induce isomorphisms on cohomology.
\end{proposition}

\begin{proof}
Consider the Laplace operators 
$\Delta_{\partial} = \partial\partial^* + \partial^*\partial$
and $\Delta_{\bar{\partial}} = \bar{\partial}\bar{\partial}^*
+ \bar{\partial}^*\bar{\partial}$.
On K\"{a}hler manifold (see e.g.
Griffiths-Harris \cite{Gri-Har}, p. 115),
we have $\Delta_{\partial} = \Delta_{\bar{\partial}} =
\frac{1}{2}  \Delta_d$.
By Hodge theory, every element in $H(\Omega^{*, *}, \bar{\partial})$
is represented by an element $\alpha$, such that
$\Delta_{\bar{\partial}}\alpha = 0$.
It then follows that $\Delta_{\partial} \alpha = 0$, 
and hence $\partial^* \alpha = 0$, 
i.e., $\alpha \in \Ker \partial^*$.
This proves that $H(i)$ is surjective.

To show that $H(i)$ is injective,
let $\alpha \in \Ker \partial^*$, 
such that $\alpha = \bar{\partial} \beta$ for some $\beta$.
We need to show that $\alpha \in \bar{\partial} \Ker \partial^*$.
By Hodge decomposition,
$\bar{\partial} \beta = \Delta_{\bar{\partial}} \gamma$ for some
$\gamma$.
Then we have
$$\Delta_{\bar{\partial}}\bar{\partial}\gamma 
= \bar{\partial}\Delta_{\bar{\partial}} \gamma 
= \bar{\partial}^2 \beta = 0,
$$
which implies that $\bar{\partial}\gamma = 0$.
Furthermore,
$$\Delta_{\partial}\partial^*\gamma = 
\partial^* \Delta_{\partial} \gamma 
= \partial^* \Delta_{\bar{\partial}} \gamma 
= \partial^* \alpha = 0$$
implies that $\partial^* \gamma = 0$.
Now $\alpha = \Delta_{\bar{\partial}}\gamma =
\bar{\partial}(\bar{\partial}^*\gamma) + 
\bar{\partial}^*(\bar{\partial}\gamma)
= \bar{\partial}(\bar{\partial}^*\gamma)$.
But
$$\partial^*(\bar{\partial}^* \gamma) 
= - \bar{\partial}^*(\partial^*  \gamma) = 0$$
show that $\bar{\partial}^* \gamma \in \Ker \partial^*$,
therefore, $\alpha \in \bar{\partial} \Ker \partial^*$,
thus  $H(i)$ is injective. 
The proof for $H(j)$ is similar.
\end{proof}

Now we use Theorem \ref{thm:construction} to get

\begin{theorem}
For any closed K\"{a}hler manifold $X$, if $K$ is the
algebra of formal power series generated by the dual space
of the Dolbeault cohomology $H^*_{\bar{\partial}}(X)$,
then there is a structure of formal Frobenius manifold 
on the formal spectrum of $K$ obtained from the dGBV
algebra $(\Omega^{*, *}(X), \wedge, \delta = \bar{\partial},
\Delta = \partial^*)$.
\end{theorem}

\begin{remark}\label{convergence}
The construction involves the solution 
by power series method to the equation
$$\bar{\partial} \Gamma + \frac{1}{2}[\Gamma, \Gamma] = 0,$$
where $\Gamma = \sum_n \Gamma_n$,
$\Gamma_0 = 0$, $\Gamma_1 = \sum x^j \eta_j$,
$\Gamma_n $ is  a (super)polynomial in $x^j$'s with coefficients in 
$\Img \bar{\partial}\partial^*$.
Here $\eta_j$ is a basis of $\Ker \bar{\partial} \cap \Ker \bar{\partial}^*
\cong H = H(\Gamma(X,\Lambda^*TX \otimes \Lambda^* \overline{T}^*X), 
\bar{\partial})$,
$x^j$ is the dual basis of $H'$.
Notice that $\Gamma$ is a (super)polynomial for  $x^j$'s 
(Fermionic coordinates) with odd degrees,
and is a power series only for $x^j$'s with even degrees (Bosonic coordinates). 
A modification of the standard argument in Kodaira-Spencer-Kuranishi deformation
theory (see e.g. Morrow-Kodaira \cite{Mor-Kod}, $\S 4.2$)
shows that these series converges for all Fermionic $x^j$'s,
and for small values of Bosonic $x^j$'s 
(this method was also used by Tian \cite{Tia1} and Todorov \cite{Tod}). 
Therefore, one gets a Frobenius supermanifold with the Bosonic  part an
open neighborhood of the zero vector in $H^{even}$. 
\end{remark}

\begin{remark}
The same construction carries through if we take
$\delta = \partial$, $\Delta  = \bar{\partial}^*$.
\end{remark}

\section{Comparison with Barannikov-Kontsevich's Frobenius manifold}

Assume that $X$ is a closed Calabi-Yau $n$-manifold,
fix a nowhere vanishing holomorphic $n$-form $\Omega \in 
\Gamma(X, \Lambda^n T^*X)$.
Barannikov-Kontsevich \cite{Bar-Kon} considered
the dGBV algebra with
\begin{eqnarray*}
& {\cal A} = \oplus {\cal A}_k, & 
{\cal A}_k = \oplus_{q+p= k} 
  \Gamma(X,\Lambda^pTX \otimes \Lambda^q \overline{T}^*X).
\end{eqnarray*}
Here  ${\cal A}$ is ${\Bbb Z}$-graded. It also has  an induced
${\Bbb Z}_2$-grading.
This grading is different from that in \cite{Bar-Kon} 
(shifted by $1$),
we adopt this grading to be compatible with the notations of 
Manin \cite{Man}. 
The multiplication $\wedge$ on ${\cal A}$ is given by
the ordinary wedge products on $\Lambda^*\overline{T}^*X$
and $\Lambda^*TX$.
The derivation is $\delta = \bar{\partial}$.
Notice that for any two integers $0 \leq p, q \leq n$,
$\Omega$ defines an isomorphism 
\begin{eqnarray*}
\Lambda^pTX \otimes \Lambda^q \overline{T}^*X 
\rightarrow 
\Lambda^{n-p}T^*X \otimes \Lambda^q \overline{T}^*X, &
\gamma \mapsto \gamma \vdash \Omega,
\end{eqnarray*}
defined by contraction of the $p$-vector with $\Omega$
to get a form of type $(n-p, 0)$.
Then $\Delta$ is defined by 
$$(\Delta \gamma) \vdash \Omega = \partial(\gamma \vdash \Omega).$$
Tian-Todorov lemma shows that the bracket,  defined by
$$[\alpha \bullet \beta]
= (-1)^{|\alpha|}(\Delta(\alpha \wedge \beta) 
- \Delta\alpha \wedge \beta - (-1)^{\alpha|}
\alpha \wedge \Delta\beta), $$
is given by the wedge product on  type $(0, *)$-forms and 
the Schouten-Nijenhuis bracket on  type $(*, 0)$ 
polyvector fields. 
It is straightforward to check that 
$({\cal A}, \wedge, \delta, \Delta)$ as above is a dGBV algebra.
Furthermore, the linear functional
\begin{eqnarray*}
\int: {\cal A} \rightarrow {\Bbb C}, & 
\int \gamma = \int_X (\gamma \vdash \Omega) \wedge \Omega
\end{eqnarray*} 
is a nice integral on $X$.

Since the complex $({\cal A}, \delta)$ contains the
deformation complex 
$$\Omega^{0, 0}(TX) \stackrel{\bar{\partial}}{\rightarrow}
\Omega^{0, 1}(TX) \stackrel{\bar{\partial}}{\rightarrow}
\Omega^{0, 2}(TX),
$$
It is called the extended deformation complex. 
Barannikov-Kontsevich \cite{Bar-Kon} 
(Lemma 2.1) showed that 
the deformation functor associated with the graded differential Lie algebra 
$({\cal A}[-1], \delta, [\cdot\bullet\cdot])$
is represented by the formal spectrum of formal power series 
generated by $H'$, the dual of $H=H({\cal A}, \delta)$.
This is the extended moduli space of complex structures in Witten \cite{Wit}.
Furthermore, they constructed a structure of a formal 
Frobenius manifold on it.

\begin{remark} 
As in Remark \ref{convergence},  One  gets Frobenius supermanifold 
structure.
\end{remark}

To facilitate the comparison with our construction,
we now rewrite the operator $\Delta$ in 
Barannikov-Kontsevich \cite{Bar-Kon},  
which appeared in earlier work of Tian \cite{Tia} and Todorov \cite{Tod}
with different notation, in terms of more familiar operator.
To start with, notice that the Hermitian metric induces
an isomorphism $TX \cong \overline{T}^*X$.
Therefore, we get, for any two integers $0 \leq p, q \leq n$,
an isomorphism 
$$f_{p, q}: \Lambda^pTX \otimes \Lambda^q \overline{T}^*X
\rightarrow  \Lambda^p \overline{T}^*X \otimes  \Lambda^qTX.$$
On $\Lambda^p \overline{T}^*X \otimes  \Lambda^qTX$,
we consider the opeartor
$$\bar{\partial}^*: 
\Gamma(X, \Lambda^p \overline{T}^*X \otimes  \Lambda^qTX)
\rightarrow 
\Gamma(X, \Lambda^{p-1} \overline{T}^*X \otimes  \Lambda^qTX).$$

\begin{lemma} 
We have the following commutative diagram:
$$\CD
\Gamma(X, \Lambda^pTX \otimes \Lambda^q \overline{T}^*X) 
@>{\Delta}>>
\Gamma(X, \Lambda^{p-1}TX \otimes \Lambda^q \overline{T}^*X) \\
@V{f_{p,q}}VV @VV{f_{p-1, q}}V \\
\Gamma(X, \Lambda^p \overline{T}^*X \otimes  \Lambda^qTX)
@>>{-\bar{\partial}^*}> 
\Gamma(X, \Lambda^{p-1} \overline{T}^*X \otimes  \Lambda^qTX).
\endCD $$
\end{lemma}

\begin{proof}
It is well-known that both $\omega$ and $\Omega$ 
are parallel with respect to the Levi-Civita connection.
It follows that near  each $x \in X$, one can find 
a local frame $\{e^1, \cdots, e^n \}$ of $T^{1, 0}X$, 
such that $\Omega = e^1 \wedge \cdots e^n$,
$\omega = e^1 \wedge \bar{e}^1 + \cdots + e^n \wedge \bar{e}^n$.
Furthermore, $\nabla_{e_j} e_k = \nabla_{\bar{e}_j}e_k 
= \nabla_{e_j} \bar{e}_k = 
\nabla_{\bar{e}_j} \bar{e}_k = 0$ at $x$.
Let $\{e^1, \cdots, e^n\}$ be the dual basis.
Near $x$, every $\gamma \in 
\Gamma(X, \Lambda^{p-1}TX \otimes \Lambda^q \overline{T}^*X)$
can be written as
$$\gamma = \sum \gamma_{JK} e_J \otimes \bar{e}^K,$$
wehre $J$ runs over multiple indices 
$(j_1, \cdots, j_p)$, $1 \leq j_1 < \cdots < j_p \leq n$,
$K$ run  over multiple inddices $(k_1, \cdots, k_q)$,
$1 \leq k_1 < \cdots < k_q \leq n$,
$e_J = e_{j_1} \wedge \cdots \wedge e_{j_p}$,
$\bar{e}^K =\bar{e}^{k_1} \wedge \cdots \wedge \bar{e}^{k_q}$.
Then we have
\begin{eqnarray*}
f_{p-1,q}^{-1} \bar{\partial}^*f_{p, q}\gamma &=&
f_{p-1,q}^{-1} (\sum \bar{\partial}^*(\gamma_{JK} \bar{e}^J 
	 \wedge e_K)) \\
& = & -f_{p-1,q}^{-1} 
(\sum (\nabla_{e_j} \gamma_{JK}) (\bar{e}_j \vdash 
	\bar{e}^J) \wedge e_K) \\
& = & -\sum (\nabla_{e_j} \gamma_{JK}) (e^j \vdash 
	e_J) \wedge \bar{e}^K.
\end{eqnarray*}
On the other hand,
\begin{eqnarray*}
(\Delta \gamma) \vdash \Omega & = & 
\partial (\gamma \vdash \Omega) \\
& = & \sum \partial(\gamma_{JK} (e_J \vdash \Omega) 
	 \wedge \bar{e}^K) \\
& = & \sum (\nabla_{e_j} \gamma_{JK}) (e^j \wedge 
	(e_J \vdash \Omega)) \wedge \bar{e}^K.
\end{eqnarray*}
It is then obvious that
$f_{p-1,q}^{-1} \bar{\partial}^*f_{p, q}\gamma
= -\Delta \gamma$.
\end{proof}

Denote $\widehat{\bar{\partial}^*} = 
f_{p-1,q}^{-1} \bar{\partial}^*f_{p, q}$.
We have the following interesting comparison of the two dGBV algebras:
\begin{align*}
&{\cal A} = \oplus_{p, q} 
\Gamma(X, \Lambda^pTX \otimes \Lambda^q\overline{T}^*X), && 
\delta = \bar{\partial}, & \Delta = - \widehat{\bar{\partial}^*}, \tag{BK}  \\
&{\cal A} = \oplus_{p, q} 
\Gamma(X, \Lambda^pT^*X \otimes \Lambda^q\overline{T}^*X), && 
\delta = \bar{\partial}, & \Delta = \partial^*, \tag{CZ}
\end{align*}
(The definitions of the multiplications are similar!) 
It is reasonable to conjecture that 
the two corresponding Frobenius manifolds are mirror images of each other.


\begin{thebibliography}{99}



\bibitem{Bar-Kon} S. Barannikov, M. Kontsevich,
{\em Frobenius Manifolds and Formality of 
Lie Algebras of Polyvector Fields},
preprint, alg-geom/9710032.


\bibitem{Ber-Pet} J. Bertin, C. Peters, 
{\em Variations de structures de Hodge, vari\'{e}t\'{e}s 
de Calabi-Yau et sym\'{e}trie miroir},
in {\bf Introduction \`{a} la th\'{e}orie de Hodge},
169--256, Panor. Synth\`{e}ses, 3, Soc. Math. France, Paris, 1996.

\bibitem{Bry} J.-L. Brylinski, 
{\em A differential complex for Poisson manifolds},
{\bf J. Differential Geom. 28} (1988), no. 1, 93--114.

 \bibitem{Can-Oss-Gre-Park} 
P. Candelas, X. C. de la Ossa, p. S. Green,  L. Parkes, 
{\em A pair of Calabi-Yau manifolds as
an exactly soluble superconformal theory}, 
{\bf Nuclear Phys. B 359} (1991), no. 1, 21--74.

\bibitem{Cao-Zho1} H.-D. Cao, J. Zhou,
{\em On Quantum de Rham Cohomology Theory}, 
preprint, April, 1998.

\bibitem{Cao-Zho2} H.-D. Cao, J. Zhou,
{\em Identification of Two Frobenius Manifolds In Mirror Symmetry}, 
preprint,  1998.

\bibitem{Dub1} B. Dubrovin,
{\em Integrable systems in topological field theory},
{\bf  Nuclear Phys. B 379} (1992), no. 3,
627--689. 

\bibitem{Dub2} B. Dubrovin,
{\em Geometry of $2$D topological field theories},
in {\bf Integrable systems and quantum groups}
(Montecatini Terme, 1993), 120--348, 
Lecture Notes in Math., 1620, Springer, Berlin, 1996. 


\bibitem{Gol-Mil} W.M. Goldman, J.J. Millson,  
{\em The deformation theory of representations of fundamental groups of
compact Kähler manifolds}, 
{\bf Inst. Hautes Études Sci. Publ. Math. No. 67} (1988), 43--96.

\bibitem{Gri-Har} P. Griffiths, J. Harris, 
{\em Principles of algebraic geometry}. 
Pure and Applied Mathematics.
Wiley-Interscience [John Wiley \& Sons], New York, 1978.

\bibitem{Kon} M. Kontsevich, 
{\em Homological algebra of mirror symmetry},
in {\bf  Proceedings of the International Congress of
Mathematicians}, 
Vol. 1, 2 (Zürich, 1994), 120--139, Birkhäuser, Basel, 1995.

\bibitem{Kon2} M. Kontsevich, 
{\em Deformation quantization of Poisson manifolds, I},
preprint,  q-alg/9709040.



\bibitem{Kos} J.-L. Koszul, 
{\em Crochet de Schouten-Nijenhuis et cohomologie},
{\bf The mathematical heritage of \'{E}lie Cartan (Lyon, 1984). 
Astérisque} 1985, Numero Hors Serie, 257--271.
 
\bibitem{Li-Tia} J. Li, G. Tian,
{\em Comparison of the algebraic and 
the symplectic Gromov-Witten invariants},
preprint, alg-geom/9712035.


\bibitem{Lia-Liu-Yau} B. Lian, K. Liu, S.T. Yau,
{\em Mirror principle I}, preprint, alg-geom/9712011.


\bibitem{Liu} G. Liu, 
{\em Associativity of quantum  multiplication}, preprint, 1994.


\bibitem{McD-Sal}
D. McDuff, D. Salamon, 
{\bf $J$-holomorphic curves and quantum cohomology}.
 University Lecture Series, 6. 
 American Mathematical Society, Providence, RI, 1994.


\bibitem{Man} Y. Manin, 
{\em Three constructions of Frobenius manifolds: a comparative study},
preprint, math.QA/9801006.


\bibitem{Mor} D.R. Morrison, 
{\em Mirror symmetry and rational curves on quintic threefolds: 
a guide for mathematicians},
{\bf J. Amer. Math. Soc. 6} (1993), no. 1, 223--247.

\bibitem{Mor-Kod} J. Morrow, K. Kodaira 
{\bf Complex manifolds},
 Holt, Rinehart and Winston, Inc., New
York-Montreal, Que.-London, 1971.

\bibitem{Rua-Tia} Y. Ruan, G. Tian, 
{\em A mathematical theory of quantum cohomology},
{\bf  J. Differential Geom. 42}
(1995), no. 2, 259--367.

\bibitem{Tia1} G. Tian,
{\em Smoothness of the universal deformation space of compact
Calabi-Yau manifolds and its Petersson-Weil metric},
in {\bf  Mathematical aspects of string theory}
(San Diego, Calif., 1986), 629--646, Adv. Ser. Math. Phys., 1, World Sci. Publishing,
Singapore, 1987.

\bibitem{Tia} G. Tian,
{\em Quantum cohomology and its associativity}, 
in {\bf  Current developments in mathematics}, 1995
(Cambridge, MA), 361--401, Internat. Press, Cambridge, MA, 1994.

\bibitem{Tod} A.N. Todorov, 
{\em The Weil-Petersson geometry of 
the moduli space of ${\rm SU}(n\geq 3)$
(Calabi-Yau) manifolds. I.},
{\bf Comm. Math. Phys. 126} (1989), no. 2, 325--346. 

\bibitem{Vaf} C. Vafa, 
{\em Topological mirrors and quantum rings},
in {\bf  Essays on mirror manifolds}, S.T. Yau ed., 
96--119, Internat. Press, Hong Kong, 1992.

\bibitem{Vai} I. Vaisman,
{\bf Lectures on the geometry of Poisson manifolds}. 
Progress in Mathematics, 118. Birkhäuser
Verlag, Basel, 1994. 

\bibitem{Wit1} E. Witten,
{\em Two-dimensional gravity and intersection theory on moduli
space}, 
in {\bf Surveys in differential geometry} (Cambridge, MA, 1990), 243--310, Lehigh Univ.,
Bethlehem, PA, 1991.

\bibitem{Wit} E. Witten, 
{\em Mirror manifolds and topological field theory},
in  {\bf Essays on mirror manifolds},
S. T. Yau ed., International Press Co., Hong Kong, 1992.

\bibitem{Yau} S.T. Yau ed., 
{\bf Essays on mirror manifolds}, 
 International Press, Hong Kong, 1992.

\end{thebibliography}
\end{document}